\documentclass[]{article}
\begin{document}
\newtheorem{proposition}{Proposition}[section]
\newtheorem{definition}{Definition}[section]
\newtheorem{lemma}{Lemma}[section]

\title{\bf Noncommutative Partial Derivative}
\author{Keqin Liu\footnote{Email address: kliu at math.ubc.ca} \\Department of Mathematics\\The University of British Columbia\\Vancouver, BC\\
Canada, V6T 1Z2}
\date{May, 2022}
\maketitle

\begin{abstract} We introduce the axiomatic definition   of the point-derivative for noncommutative algebras  and present the counterparts of  the ordinary multi-variable chain rule and Clairaut's Theorem in the context of partial point-derivatives.

\end{abstract}

\bigskip
Among various attempts of extending the notion of the ordinary derivative to noncommutative polynomials, there are two noteworthy extensions of the  ordinary derivative. One is the Hausdorff derivative, and the other one is the cyclic derivative introduced in \cite{RSS}. Some well-known properties of  the ordinary derivative can be naturally extended to 
the cyclic derivative for formal series in  a single variable,   but any significant properties of 
cyclic derivative for formal series in several variables  has not been obtained  since this problem was  mentioned in \cite{RSS}. In particular, there is not any  counterpart of  the ordinary multi-variable chain rule and Clairaut's Theorem in the context of both Hausdorff derivatives and   cyclic derivatives.

\medskip
In this paper, we introduce
the axiomatic definition of  the point-derivative for  noncommutative algebras,  explain how to define 
the  point-derivative  in the algebra of noncommutative formal power series, and  present the 
counterparts of  the ordinary multi-variable chain rule and Clairaut's Theorem in the context of partial point-derivatives.

\medskip
Throughout this paper, we use the same notations and terminologies as the ones in  \cite{RSS}.  In particular, $K$ denotes a field of 
characteristic zero and all algebras are over the field $K$ and have the identity.  

\bigskip
\section{The Notion of  Point-derivative}

Let $\cal{A}$ be a $K$-algebra. Recall that a $K$-linear map $d: \cal{A}\to\cal{A}$  is called  a {\it derivation} on $\cal{A}$ if $d(ab)=a d(b)+d(a) b$ for $a$, $b\in \cal{A}$. The axiomatic definition of  the point-derivative for  noncommutative algebras  is given in the following 

\medskip
\begin{definition}\label{def1.1} Let $\cal{A}$ be a $K$-algebra   and let $End(\cal{A})$ be the  set of all $K$-linear maps from $\cal{A}$ to $\cal{A}$. A $K$-linear map $\cal{D}: \cal{A}\to$  $End({\cal A})$  is called a 
{\bf point-derivation} on $\cal{A}$ if 
 \begin{description} 
\item[(i)] $\cal{D}_{\beta}:=\cal{D}(\beta)\in$ $End(\cal{A})$ is a derivation on $\cal{A}$ for all $\beta\in \cal{A}$,
\item[(ii)] ${\cal D}_z=z\,{\cal D}_1$ for all $z$ in the center of ${\cal A}$,
\end{description}
where $1$ is the identity of the algebra ${\cal A}$. For $f\in {\cal A}$ and $\beta\in \cal{A}$, 
${\cal D}_{\beta}(f)\in  {\cal A}$ is called the {\bf point-derivative depending on $\beta$} of $f$ or  {\bf $\beta$-derivative} of $f$  .
\end{definition}

\medskip
If  $\cal{A}$ is a commutative algebra, then  a  derivation $d$ on ${\cal A}$ produces naturally a  point-derivation ${\cal D}$ on 
${\cal A}$ as follows:
$$ {\cal D}(\beta):=\beta d\quad\mbox{for $\beta\in {\cal A}$},$$
where $\beta d\in End({\cal A})$ is defined by $(\beta d)(a): =\beta\cdot d(a)$ for $a\in {\cal A}$.

\bigskip
We now explain how to introduce  point-derivative in noncommutative  formal power series. Let $A$ be the alphabet set which has a distinguished letter, denoted by $x$ and called the {\it variable},  and an infinite supply of other letters $a$, $b$, \dots called 
{\it constants}.  A {\it word} $w$ is an element of the free monoid  $M$ generated by $A$, and  the {\it degree}
of a word $w$  is the number of occurrences of the letter $x$ in $w$.
The {\it empty word}, which is denoted by $1$, is the identity element of free monoid $M$. Following \cite{RSS}, we  use   
$K\{\{a, b, \dots, x\}\}$ to denote the algebra of noncommutative formal power series in a single 
variable $x$ over the field $K$.  

\medskip
Let $\beta\in K\{\{a, b, \dots, x\}\}$ be a fixed formal power series.  For a given word 
$w=c_1x^{i_1}c_2x^{i_2}\cdots c_nx^{i_n}c_{n+1}\in M$, where $c_1$, $c_2$, $\cdots$, $c_n$, $c_{n+1}$ are the words of degree $0$ and $i_1$, $i_2$, $\cdots$, $i_n$ are nonnegative integers, we define a map  $\displaystyle\frac{d}{d_{\beta}x}$ from $M$ to $K\{\{a, b, \dots, x\}\}$ as follows:
$$
\displaystyle\frac{d}{d_{\beta}x}(w):=0\quad\mbox{if $i_1=i_2=\cdots =i_n=0$}
$$
and
\begin{eqnarray*}
\displaystyle\frac{d}{d_{\beta}x}(w):&=&\displaystyle\sum_{k=1}^nc_1x^{i_1}\cdots c_{k-1}x^{i_{k-1}}
c_k\cdot\displaystyle\frac{d}{d_{\beta}x}(x^{i_k})\cdot  c_{k+1}x^{i_{k+1}}\cdots c_nx^{i_n}c_{n+1}\\
&&\mbox{if $i_1i_2\cdots i_n> 0$,}
\end{eqnarray*}
where $c_1x^{i_1}\cdots c_{k-1}x^{i_{k-1}}:=1$ for $k=1$,  $c_{k+1}x^{i_{k+1}}\cdots c_nx^{i_n}c_{n+1}:=1$ for $k=n$ and 
 $\displaystyle\frac{d}{d_{\beta}x}(x^{i_k})$ is defined  by
$$
\displaystyle\frac{d}{d_{\beta}x}(x^{i_k})=\left\{\begin{array}{ll}\beta&
\mbox{if $i_k=1$}\\
x^{i_k-1}\beta+x^{i_k-2}\beta x+\cdots +x\beta x^{i_k-2}+\beta x^{i_k-1}&\mbox{if $i_k> 1$}\end{array}\right..
$$
After extending the map $\displaystyle\frac{d}{d_{\beta}x}:M\to K\{\{a, b, \dots, x\}\}$ linearly and continuously, we get a 
$K$-linear map $\displaystyle\frac{d}{d_{\beta}x}: K\{\{a, b, \dots, x\}\}\to K\{\{a, b, \dots, x\}\}$, which is called the
{\bf point-derivative operator depending on $\beta$}. Clearly, The map: 
$\beta\to \displaystyle\frac{d}{d_{\beta}x}$ is a point-derivation on $K\{\{a, b, \dots, x\}\}$.   
The  {\bf $\beta$-derivative}  $\displaystyle\frac{d}{d_{\beta}x}\big(f(x)\big)$
  of a formal power series $f(x)\in K\{\{a, b, \dots, x\}\}$ is also denoted by  
$\displaystyle\frac{df}{d_{\beta}x}$  or $f'_{\beta}(x)$.
 Note that  the $1$--derivative $\displaystyle\frac{df}{d_{1}x}=f'_{1}(x)$ of the  formal power series $f(x)$ is the Hausdorff derivative $H<f>$ of $f(x)$.

\medskip
Like the Hausdorff derivative and cyclic derivative, we have the following characterization of the  point-derivative operator 
 $\displaystyle\frac{d}{d_{\beta}x}$.

\begin{proposition}\label{pr1.1} Let  $\beta\in K\{\{a, b, \dots, x\}\}$ be a fixed formal power series. The point--derivative operator 
 $\displaystyle\frac{d}{d_{\beta}x}:  K\{\{a, b, \dots, x\}\}\to K\{\{a, b, \dots, x\}\}$ is the unique $K$-linear continuous operator on 
$ K\{\{a, b, \dots, x\}\}$ satisfying 
\begin{description} 
\item[(i)] $\displaystyle\frac{d}{d_{\beta}x}(a)=0$ for any constant $a$,
\item[(ii)] $\displaystyle\frac{d}{d_{\beta}x}(x)=\beta$,
\item[(iii)] $\displaystyle\frac{d}{d_{\beta}x}(fg)=\displaystyle\frac{d}{d_{\beta}x}(f)\cdot g+f\cdot\displaystyle\frac{d}{d_{\beta}x}(g)$ for  $f$,  $g\in  K\{\{a, b, \dots, x\}\}$.
\end{description}
\end{proposition}

\bigskip
\section{Partial Point-derivatives}

For convenience, we use  $K\{\{a, b, \dots, || x, y\}\}$ to denote the algebra of noncommutative formal power series in two  noncommutative  variables $x$ and $y$.  Let   $f(x, y)\in K\{\{a, b, \dots, || x, y\}\}$  and 
$\beta=\beta(x, y)\in K\{\{a, b, \dots, || x, y\}\}$ be   formal power series in two  noncommutative  variables $x$ and $y$. 
 If we keep $y$ constant, then  $f(x, y)$ and $\beta=\beta(x, y)$ are  formal power series of $x$ and its  $\beta$-derivative  with respect to the variable $x$ is denoted by $ \displaystyle\frac{\partial f(x, y)}{\partial_{\beta} x}$ or  $ \displaystyle\frac{\partial f}{\partial_{\beta} x}$, which is called 
the {\bf partial  $\beta$-derivative of $f(x, y)$  with respect to the variable $x$}.  
 If we keep $x$ constant, then $f(x, y)$  and $\beta=\beta(x, y)$ are  formal power series of $y$ and its  $\beta$-derivative  with respect to the variable $y$ is denoted by $ \displaystyle\frac{\partial f(x, y)}{\partial_{\beta} y}$ or  $ \displaystyle\frac{\partial f}{\partial_{\beta} y}$, which is called 
the {\bf partial  $\beta$-derivative of $f(x, y)$  with respect to the variable $y$}.  

\medskip
Let $\beta=\beta(x, y)$, $\gamma=\gamma (x, y)$ and  $f=f(x, y)$ be  formal 
power series in   noncommutative  variables $x$ and $y$.    We define  {\bf the second partial 
 $( _{\gamma}y,  \, _{\beta}x)$-derivative 
$\displaystyle\frac{\partial^2 f}{\partial_{\gamma} y\partial_{\beta} x}$ of $f(x, y)$} and
 {\bf the second partial 
 $( _{\beta}x,  \, _{\gamma}y)$-derivative 
$\displaystyle\frac{\partial^2 f}{\partial_{\beta} x\partial_{\gamma} y}$ of $f(x, y)$} by
$$
\displaystyle\frac{\partial^2 f}{\partial_{\gamma} y\partial_{\beta} x}:=
\displaystyle\frac{\partial }{\partial_{\gamma} y}\left(\displaystyle\frac{\partial f}{\partial_{\beta} x}\right)
\quad\mbox{and}\quad
\displaystyle\frac{\partial^2 f}{\partial_{\beta} x\partial_{\gamma} y}:=
\displaystyle\frac{\partial }{\partial_{\beta} x}\left(\displaystyle\frac{\partial f}{\partial_{\gamma} y}\right).
$$

\medskip
The following proposition extends  multi-variable chain rule 
and Clairaut's Theorem to  partial point-derivatives for formal power series  containing $x$ and $y$.  

\medskip
\begin{proposition}\label{pr2.3} If $f=f(x, y)$,  $\beta=\beta(x, y)$, $\gamma=\gamma (x, y)$, 
$u=u(x, y)$ and $v=v(x, y)$ are  formal power series in two  noncommutative  variables $x$ and $y$, then
$$
\displaystyle\frac{\partial \big(f(u, v)\big)}{\partial_{\beta} x}=
\displaystyle\frac{\partial \big(f(u, v)\big)}{\partial_{\frac{\partial u}{\partial_{\beta} x}} u}+
\displaystyle\frac{\partial \big(f(u, v)\big)}{\partial_{\frac{\partial v}{\partial_{\beta} x}} v}.
$$
and
$$
\displaystyle\frac{\partial^2 f}{\partial_{\gamma} y\partial_{\beta} x}-
\displaystyle\frac{\partial^2 f}{\partial_{\beta} x\partial_{\gamma} y}=
\displaystyle\frac{\partial f}{\partial_{\big(\frac{\partial \beta}{\partial_{\gamma} y}\big)} x}-
\displaystyle\frac{\partial f}{\partial_{\big(\frac{\partial \gamma}{\partial_{\beta} x}\big)} y}.
$$
\end{proposition}

\medskip
The proof of the proposition above follows from a direct computation.

\bigskip
By \cite{CKW}, the ordinary partial derivatives is a useful tool in the study of commutative arithmetic
circuits. Hence, one problem, which is worthy to study it in future, is to determine if  the
partial point-derivatives can be used to study non-commutative arithmetic circuits.


\bigskip
\begin{thebibliography}{9}
\bibitem{CKW} Xi Chen, Neeraj Kayal and Avi Wigderson, \textsl{Partial Derivatives in Arithmetic Complexity and Beyond}, Foundations and Trends in Theoretical Computer Science, 
Vol.6, Nos.1-2(2010) 1-138.
\bibitem{RSS} Gian-Carlo Rota,  Bruce Sagan \& Paul R. Stein, \textsl{A Cyclic Derivative in Noncommutative Algebra}, 
Journal of Algebra 64, 54-75 (1980).
\end{thebibliography}
\end{document}